\documentclass[11pt,reqno]{amsart}
\usepackage{enumerate}
\usepackage{amsmath, amsfonts, amssymb, amsthm}
\usepackage{verbatim}
\usepackage[mathscr]{eucal}

\def\l{\left}
\def\r{\right}
\def\bg{\bigg}
\def\({\bg(}
\def\){\bg)}
\def\[{\bg[}
\def\]{\bg]}
\def\t{\mbox}
\def\f{\frac}

\def\em{\emptyset}
\def\se {\subseteq}

\def\sm{\setminus}

\def\bi{\binom}
\def\eq{\equiv}

\def\mo{\t{\rm mod}}

\def\ord{\t{\rm ord}}

\def\exp{\t{\rm exp}}

\def\sss{{\mathsf s}}

\newcommand{\N}{\mathbb N}
\newcommand{\Z}{\mathbb Z}

\newcommand{\Q}{\mathbb Q}

\newcommand{\Fc}{\mathcal F}
\newcommand{\vp}{\mathsf v}

\newcommand{\supp}{\text{\rm supp}}

\numberwithin{equation}{section}
\newtheorem{Theorem} {Theorem} [section]

\newtheorem{Lemma} [Theorem] {Lemma}

\theoremstyle{definition}

\newcommand{\Sum}[2]{\sum_{#1}^{#2}}
\newcommand{\Summ}[1]{\sum_{#1}}

 \newcommand{\sbinom}[2]{\genfrac{[}{]}{0pt}{}{#1}{#2}}

\newcommand{\be}{\begin{equation}}
\newcommand{\ee}{\end{equation}}
\newcommand{\ber}{\begin{eqnarray}}
\newcommand{\eer}{\end{eqnarray}}
\newcommand{\nn}{\nonumber}

\newcommand{\und}{\;\mbox{ and } \;}

\numberwithin{equation}{section}

\newcommand{\la}{\langle}
\newcommand{\ra}{\rangle}

\newcommand{\rk}{\mathsf{rk}\,}

\begin{document}
\hbox{Accepted by Adv. in Appl. Math.}
\medskip
\title
[On Weighted Zero-sum Sequences]{On Weighted Zero-sum Sequences}
\author[Sukumar Das Adhikari, David J. Grynkiewicz, and Zhi-Wei Sun]
{Sukumar Das Adhikari, David J. Grynkiewicz and Zhi-Wei Sun$^{1}$}

\thanks{$^{1}$Corresponding author. Supported by the National Natural
Science Foundation (grant 11171140) of China.}
\address{Harish-Chandra Research Institute, Chhatnag Road,  Jhusi,
Allahabad 211 019, India} \email{adhikari@mri.ernet.in}
\address{Institut f\"ur Mathematik und Wissenschaftliches  Rechnen.
Karl-Franzens-Universit\"at. Heinrichstrasse 36. 8010 Graz, Austria}
\email{diambri@hotmail.com}
\address{Department of Mathematics, Nanjing University,  Nanjing 210093,
People's Republic of China} \email{zwsun@nju.edu.cn}


\begin{abstract}
Let $G$  be a finite additive abelian group with exponent $\exp(G)=n>1$
and let $A$ be a nonempty subset of $\{1,\ldots,n-1\}$.
In this  paper, we investigate the smallest positive integer $m$, denoted by
$\mathsf s_A(G)$,
such that any sequence $\{c_i\}_{i=1}^m$ with terms from $G$
has a length $n=\exp(G)$ subsequence $\{c_{i_j}\}_{j=1}^{n}$  for which
there are $a_1,\ldots,a_n\in A$ such that $\sum_{j=1}^na_ic_{i_j}=0$.

When $G$  is a $p$-group, $A$ contains no multiples of $p$ and any
two distinct elements of $A$ are incongruent  mod $p$, we show that
$\mathsf s_A(G)\le \lceil \mathsf{D}(G)/|A|\rceil+\exp(G)-1$ if
$|A|$ is at least $(\mathsf{D}(G)-1)/(\exp(G)-1)$, where $\mathsf{D}(G)$ is the
Davenport constant of $G$ and this upper bound for $\mathsf s_A(G)$
in terms of $|A|$ is essentially best possible.

In the case $A=\{\pm 1\}$, we determine the asymptotic behavior of $\mathsf s_{\{\pm 1\}}(G)$ when $\exp(G)$ is even,
showing that,
for  finite abelian groups of even exponent and fixed rank,
$$\mathsf s_{\{\pm 1\}}(G)=\exp(G)+\log_2|G|+O(\log_2\log_2|G|)\;\ \ \mbox{ as }\;\exp(G)\rightarrow +\infty.$$
Combined with a lower bound of $\exp(G)+\Sum{i=1}{r}\lfloor\log_2
n_i\rfloor$, where $G\cong \Z_{n_1}\oplus\cdots\oplus \Z_{n_r}$ with
$1<n_1|\cdots |n_r$, this determines $\mathsf s_{\{\pm 1\}}(G)$, for
even exponent groups, up to a small order error term. Our method
makes use of the theory of $L$-intersecting set systems.

Some additional more specific values and results related to $\mathsf s_{\{\pm 1\}}(G)$ are also computed.
\end{abstract}

\maketitle

\section{Introduction}

Let $G$ be a finite abelian group  written additively and let $\mathscr F(G)$
be the set of all finite,  ordered  sequences with terms from  $G$, though
the ordering will not be of relevance  to our investigations apart from
notational concerns.
A sequence $S=\{c_i\}_{i=1}^n\in \mathscr F(G)$ is said to be a zero-sum
sequence if $\sigma(S):=c_1+\cdots+c_n=0$.
In the theory of zero-sums, the constant $\sss(G)$ is defined to  be the
smallest positive integer $n$ such that
any sequence of length $n$ contains a zero-sum subsequence of  length
$\exp(G)$ (the exponent of $G$).  By \cite[Theorem 6.2]{GG}, we have $\mathsf
s(G)\le|G|+\exp(G)-1$.
For $n\in\Z^+=\{1,2,3,\ldots\}$, let $\Z_n=\Z/n\Z$ denote the  ring of
residue classes modulo $n$.
The famous Erd\H{o}s-Ginzburg-Ziv  Theorem (EGZ) \cite{EGZ} (see also
\cite{GG} and \cite{Nat}) implies $\sss(\Z_n)=2n-1$,
and the Kemnitz-Reiher Theorem \cite{Rei} states that $\sss(\Z_n^2)=4n-3$
where $\Z_n^2=\Z_n\oplus\Z_n$.

Shortly after the confirmation of Caro's weighted EGZ conjecture \cite{wegz},
which introduced the idea of considering certain weighted subsequence sums, Adhikari and his collaborators
(cf. \cite{AC} \cite{ACFKP} \cite{AR})
 initiated the study of a new kind of
weighted zero-sum problem.
Let $A$ be a nonempty subset  of $[1,\exp(G)-1]=\{1,\ldots,\exp(G)-1\}$.
For a  sequence $\{c_i\}_{i=1}^n\in\mathscr F(G)$,
if there are $a_1,\ldots,a_n\in  A$ such that $\sum_{i=1}^na_ic_i=0$, then
the sequence is said to have $0$ as an $A$-weighted sum or, simply, to be an $A$-weighted zero-sum subsequence.
Similar to the classical case  with $A=\{1\}$,  various $A$-weighted
constants can be defined as follows:

\begin{itemize}
\item $\mathsf D_A(G)$ is the  least integer $n$ such that any $S\in
{\mathscr F}(G)$ of length $|S|\geq n$  contains a nonempty
$A$-weighted zero-sum subsequence.
\item $\mathsf E_A(G)$ is the  least integer $n$ such that any $S\in
{\mathscr F}(G)$ with length $|S|\geq n$ has an $A$-weighted zero-sum subsequence of length $|G|$.
\item $\mathsf s_A(G)$ is the  least integer $n$ such that any $S\in
{\mathscr F}(G)$ with length $|S|\geq n$ has an $A$-weighted zero-sum subsequence
of length $\exp(G)$.
\end{itemize}

The conjecture that $\mathsf E_A(G)=|G|+\mathsf D_A(G)-1$
was recently confirmed \cite{E=G+D-1}, rendering  the independent study of
$\mathsf D_A(G)$ and $\mathsf E_A(G)$ no longer necessary.
See also \cite{AC}, \cite{weight-Gao-cyclic} and \cite{XL} for  previous
partial results on the conjecture.

Let $n$ and $r$ be positive integers. In \cite{ABPP},  Adhikari and
his coauthors investigated $$f_A(n,r):=\mathsf s_A(\Z_n^r)$$ and
proved that $f_{\{\pm1\}}(n,2) = 2n-1$ when $n$ is odd. If  $p$ is a
prime, $A\se[1,p-1]$, and  $\{a\ \mo\ p:\ a\in A\}$ is a subgroup of
the multiplicative group $\Z_p^*=\Z_p\sm\{0\}$, then the authors in
\cite{AAS} showed that
  $$f_A(p,r) \leq \frac{r(p-1)}{|A|} + p\ \ \ \t{for}\ 1\le r
<\frac{p|A|}{p-1};$$
in particular, $f_A(p,|A|) \leq 2p-1$ for such $A$.

\medskip

In the present paper, we obtain an essentially sharp upper bound for
 $\mathsf s_A(G)$---without the restriction that $\{a\ \mo\ p:\ a\in A\}$
 forms a subgroup of $\Z_p^*$---which is valid for an arbitrary abelian
$p$-group $G$.

 For an abelian $p$-group $G\cong \Z_{p^{k_1}}\oplus\cdots\oplus
\Z_{p^{k_r}}$ with $k_1,\ldots,k_r\in\N$,
 Olson \cite{Olson} proved that the Davenport constant $\mathsf
D(G)=\mathsf D_{\{1\}}(G)$ equals $\mathsf d^*(G)+1$, where
$$\mathsf d^*(G):=\sum_{t=1}^r(p^{k_t}-1).$$

Our first main theorem is as follows.

 \begin{Theorem}\label{Gen} Let $p$ be a prime and let $G$ be an abelian
$p$-group with $|G|>1$.
  Let $\em\not=A\se[1,p^{k_r}-1]\sm p\Z$ and suppose that any two distinct
elements of $A$ are incongruent modulo $p$.
 Then, for each $k\in\Z^+$, any sequence
in ${\mathscr F}(G)$ of length at least
$p^k-1+\lceil(d^*(G)+1)/|A|\rceil$ contains a nonempty  $A$-weighted
zero-sum subsequence whose length is divisible by $p^k$.
 Thus, if $|A|\ (\exp(G)-1)\ge d^*(G)=\mathsf{D}(G)-1$ (which happens if
$|A|$ is at least $\rk(G)=r$, the rank of $G$), then we have
$$\mathsf s_A(G)\le \exp(G)-1+\bg\lceil\f{\mathsf{D}(G)}{|A|}\bg\rceil.$$

\end{Theorem}

For any abelian $p$-group $G$, our upper bound  for $\mathsf s_A(G)$ in
terms of $|A|$ is essentially best
possible, as illustrated by the following example (see also
 \cite{AAS} for the particular case $G=\Z_p^r$). Note that the condition
$|A|\ge (\mathsf{D}(G)-1)/(\exp(G)-1)$ cannot be removed even in the
classical case $A=\{1\}$  since it is known that
$\mathsf s(\Z_p^2)=4p-3>\mathsf{D}(\Z_p^2)+\exp(\Z_p^2)-1=3p-2$.

\medskip
\noindent  {\it Example}. Let $p$ be a  prime and let
$G\cong\Z_{p^{k_1}}\oplus\cdots\oplus\Z_{p^{k_r}}$,
where $1\le k_1\le\cdots \le k_r$. Set
$A=[1,l]$ with $l\le p^{k_r}-1$. Consider a sequence  $S$ over $G$
which consists of
\begin{align*}&(0,0,\ldots,0)\quad\mbox{ used }\;  p^{k_r}-1\ \mbox{times},
\\&(1,0,\ldots,0) \quad \mbox{ used }\;\l\lfloor\f{p^{k_1}-1}l\r \rfloor\
\mbox{times},
\\&(0,1,\ldots,0) \quad \mbox{ used }\;\l\lfloor\f{p^{k_2}-1}l\r \rfloor\
\mbox{times},
\\&\vdots
\\&(0,\ldots,0,1) \quad \mbox{ used }\;\l\lfloor\f{p^{k_r}-1}l\r \rfloor\
\mbox{times}.
\end{align*}
Clearly, $S$ contains no subsequence of length $\exp(G)=p^{k_r}$  which
has $0$ as an $A$-weighted sum.
Note that the length of $S$ is $p^{k_r}-1+\sum_{t=1}^r
\lfloor(p^{k_t}-1)/l\rfloor$, which coincides with
 $p^{k_r}-2+\lceil (d^*(G)+1)/l\rceil=\exp(G)-2+\lceil
 \mathsf{D}(G)/|A|\rceil$ when $l$ divides  every $p^{k_t}-1$, which may
easily be arranged, for instance, if all the $k_t$ are equal.

\medskip

Under the conditions of Theorem \ref{Gen}, Thangadurai \cite{Thanga}
showed that $\mathsf{D}_A(G)\le\lceil \mathsf{D}(G)/|A|\rceil$
via the group-ring method. This result is an easy consequence of Theorem \ref{Gen} since we may add $\exp(G)-1$
0's to a sequence in $\mathscr F(G)$ of length $\lceil \mathsf{D}(G)/|A|\rceil$ and then apply our theorem.
\medskip

As we have already mentioned, in \cite{ABPP}
it was proved that $\mathsf s_{\{\pm 1\}}(\Z_n^2) = 2n-1=2\ \exp(G)-1$ if $n$ is odd.
It is easy to see (and is a specific case in the Kemnitz-Reiher Theorem \cite{Rei}
mentioned before) that $\mathsf s_{\{\pm1\}}(\Z_2^2) =5=2\ \exp(G)+1$.

In contrast to these results, in this paper we fully determine
the asymptotic behavior of $\mathsf s_{\{\pm 1\}}(G)$ when $\exp(G)$ is even,
 showing that,  for  finite abelian groups of even exponent and fixed rank,
 \be\label{starseam}\mathsf s_{\{\pm 1\}}(G)=\exp(G)+\log_2|G|+O(\log_2\log_2|G|)\;\ \ \mbox{ as }\;
 \exp(G)\rightarrow \infty.\ee
More precisely, we establish the following theorem.

\begin{Theorem}\label{thm-even-assymptotics} Let $r\geq 1$ be an integer.
Then there exists a constant $C_r$, dependent only on $r$,
such that $$\mathsf s_{\{\pm 1\}}(G)\leq \exp(G)+\log_2|G|+C_r\log_2\log_2|G|$$ for every finite abelian group $G$
of rank $r$ and even exponent.
\end{Theorem}

 In view of the lower bound on $\mathsf s_{\{\pm 1\}}(G)$ shown below in
Theorem \ref{cor-basic-bounds-dav}, Theorem \ref{thm-even-assymptotics}
 determines the value of $\mathsf s_{\{\pm 1\}}(G)$ up to the very small
order error term given in \eqref{starseam}.
 Our method makes use of fundamental results from the theory of
$L$-intersecting set systems and could be used to explicitly estimate
 the coefficient $C_r$ in specific cases as well as give bounds for how long
a  sequence $S\in \mathscr F(G)$ must be to ensure a $\{\pm 1\}$-weighted
 zero-sum subsequence of even length $n$, where $n$ is any even integer at
 least $(r+1)2^{r+1}$ with $r=\rk(G)$. To illustrate this point, and to gently
accustom the reader to the method in a more concrete setting,
 we first calculate some specific values of $\mathsf s_{\{\pm 1\}}(G)$ for
small $|G|$, and as a by-product of this investigation,
 obtain the following bounds on the weighted Davenport constant in the case
$A=\{\pm 1\}$.
 Note that $$\lfloor\log_2 |G|\rfloor=\lfloor\log_2(n_1n_2\cdots
n_r)\rfloor=\left\lfloor\Sum{i=1}{r}\log_2 n_i\right\rfloor,$$
 so that the difference between the upper and lower bounds given below is at
most $r$. In the case of cyclic groups, i.e.,
rank $r=1$, and $2$-groups, this means that equality holds. For the cyclic
case, this was first shown in \cite{ACFKP}.
 The results obtained for small $|G|$ can be combined with an inductive
argument to yield a simpler upper bound for rank $2$ groups,
which we handle in brevity at the end of Section \ref{sec-small-values}.

\begin{Theorem}\label{cor-basic-bounds-dav}Let $G$ be  a finite abelian
group with $G\cong \Z_{n_1}\oplus \Z_{n_2}\oplus \cdots\oplus \Z_{n_r}$,
where $1<n_1|\ldots|n_r$. Then $$\Sum{i=1}{r}\lfloor\log_2
n_i\rfloor+1\leq \mathsf D_{\{\pm 1\}}(G)\leq  \lfloor\log_2 |G|\rfloor+1$$
and
$$\mathsf s_{\{\pm 1\}}(G)\geq n_r+\mathsf
 D_{\{\pm 1\}}(G)-1\geq \exp(G)+\Sum{i=1}{r}\lfloor\log_2 n_i\rfloor.$$
\end{Theorem}

The next section is devoted to our algebraic proof of Theorem 1.1.
In Section 3 we introduce some terminology and notation for later use.
Section 4 contains some results more general than Theorem 1.3. In Section 5
we prove Theorem 1.2 with the help of some deep results from extremal set theory.

\section{Proof of Theorem 1.1}

The following result is well known, see, e.g., \cite[pp.
878-879]{AS}.

\begin{Lemma} {\rm (Lagrange's interpolation formula)}
Let $P(x)$ be a polynomial over the field of complex numbers, and
let $x_1,\ldots,x_n$ be $n$ distinct complex numbers. If $\deg P<n$,
then
$$P(x)=\sum_{j=1}^nP(x_j)\prod^n_{i=1\atop i\not=j}\f{x-x_i}{x_j-x_i}.$$
\end{Lemma}

We also need the following useful lemma.
As shown in \cite{Sun09}, it is very helpful when one
wants to establish certain zero-sum results for abelian $p$-groups
without appeal to the group-ring method.
\begin{Lemma}{\rm  (\cite[Lemma 4.2]{Sun09})} Let $p$ be a prime and let
$k\in\N=\{0,1,2,\ldots\}$ and $m\in\Z$ be integers. Then
$$\bi{m-1}{p^k-1}\eq \begin{cases}1\ (\mo\ p)&\mbox{if}\ p^k\mid m,\\0\ (\mo\
p)&\mbox{otherwise}.\end{cases}$$
\end{Lemma}

For convenience, for a polynomial $f(x_1,\ldots,x_n)$ over a field,
we use
 $[x_1^{k_1}\cdots x_n^{k_n}] f(x_1,\ldots,x_n)$ to denote the coefficient of
the monomial $x_1^{k_1}\cdots x_n^{k_n}$ in $f(x_1,\ldots,x_n)$.

Recall that a rational number is an $p$-adic integer (where $p$ is a
prime) if its denominator is not divisible by $p$.

\medskip
\begin{proof}[Proof of Theorem \ref{Gen}] Suppose that $G\cong
\Z_{p^{k_1}}\oplus\cdots\oplus\Z_{p^{k_r}}$ with
$1\le k_1\le\cdots\le k_r$.  Then $\mathsf d^*(G)=\sum_{t=1}^r(p^{k_t}-1)$
and $\exp(G)=p^{k_r}$.  Let $\{c_s\}_{s=1}^n\in\mathscr F(G)$ with
 $n=p^{k}-1+\lceil (\mathsf{d}^*(G)+1)/|A|\rceil$, where $k\in\Z^+$. We may
identify each $c_s$ with a vector
 $$\langle c_{s1}\ \mo\ p^{k_1},\ldots,c_{sr}\ \mo\
p^{k_r}\rangle\in\Z_{p^{k_1}}\oplus\cdots\oplus \Z_{p^{k_r}},$$
where $c_{s1},\ldots,c_{sr}$ are suitable integers. Set
 $$P(x)=\prod_{a\in A}(x-a)\in \Z[x]\ \ \ \mbox{and}\ \ \
c=\f{(-1)^{\mathsf{d}^*(G)+p^k-1}}{P(0)^n}$$ and define
 $$f(x_1,\ldots,x_n)
 =\bi{\sum_{i=1}^nP(x_i)-nP(0)-1}{p^k-1}
 \prod_{t=1}^r\bi{\sum_{s=1}^nc_{st}x_s-1}{p^{k_t}-1}-c\prod_{i=1}^nP(x_i)\in
\mathbb{Q}[x_1,\ldots,x_n].$$
 Since
 $n|A|>\mathsf{d}^*(G)+|A|(p^k-1)$,
 we have
  $$[x_1^{|A|}\cdots x_n^{|A|}]f(x_1,\ldots,x_n)=-c\und
\deg f=n|A|.$$ As $A\cap p\Z=\emptyset$ by the hypothesis,
$P(0)\not\eq0\ (\mo\ p)$ and hence $c$ is a $p$-adic integer.

 Since $P(x)\in\Z[x]$ is monic, for each $j\in\N$, there are $q_j(x),r_j(x)\in\Z[x]$ such that
 $$x^j=xP(x)q_j(x)+r_j(x)\ \ \t{and}\ \deg r_j\le \min\{j,\deg P\}.$$
 Note that $\deg (xP(x)q_j(x))=\deg(x^j-r_j(x))\le j$. Write
$$f(x_1,\ldots,x_n)=\sum_{j_1,\ldots,j_n\ge0}f_{j_1,\ldots,j_n}\prod_{i=1}^nx_i^{j_i}\in\Q[x_1,\ldots,x_n].$$
Then, as in the proof of Alon's Combinatorial Nullstellensatz
\cite{Alon}, we have
 \begin{align*}f(x_1,\ldots,x_n)=&\sum_{j_1,\ldots,j_n\ge0}f_{j_1,\ldots,j_n}\prod_{i=1}^n\l(x_iP(x_i)q_{j_i}(x_i)+r_{j_i}(x_i)\r)
\\=&\sum_{i=1}^nx_iP(x_i)h_i(x_1,\ldots,x_n)+\bar f(x_1,\ldots,x_n),
\end{align*}
where $h_i(x_1,\ldots,x_n)\in \Q[x_1,\ldots,x_n]$ (the $h_i$ here can be variously chosen), $\deg h_i+\deg
(x_iP_i(x))\le \deg f$, and $\bar f$ is given by
$$\bar f(x_1,\ldots,x_n)=\sum_{j_1,\ldots,j_n\ge0}f_{j_1,\ldots,j_n}\prod_{i=1}^nr_{j_i}(x_i).$$
Clearly,
$$f(a_1,\ldots,a_n)=\bar f(a_1,\ldots,a_n)\ \quad\mbox{for all}\ a_1,\ldots,a_n\in A',$$
where $A'=A\cup\{0\}$. Recall $\deg f=n|A|$. Thus, as $\deg h_i+\deg
(x_iP_i(x))\le \deg f$, it follows that
$$[x_1^{|A|}\cdots x_n^{|A|}]x_iP(x_i)h_i(x_1,\ldots,x_n)=[x_1^{|A|}\cdots x_n^{|A|}]x_i^{|A|+1}h_i(x_1,\ldots,x_n)=0,$$
whence
$$[x_1^{|A|}\cdots x_n^{|A|}]\bar f(x_1,\ldots,x_n)=[x_1^{|A|}\cdots x_n^{|A|}]f(x_1,\ldots,x_n)=-c.$$

As $\deg r_{j}(x)\le\deg P(x)$ for all $j\in\N$, the degree of $\bar
f$ in $x_i$ does not exceed $\deg P=|A'|-1$ for any $i=1,\ldots,n$.
Applying Lagrange's interpolation formula $n$ times, we obtain
\begin{align*}\bar f(x_1,\ldots,x_n)=&\sum_{a_n\in A'}\bar f(x_1,\ldots,x_{n-1},a_n)\prod_{b\in A'\sm\{a_n\}}\f{x_n-b}{a_n-b}
\\=&\cdots=\sum_{a_1,\ldots,a_n\in A'}\bar f(a_1,\ldots,a_n)\prod_{j=1}^n\prod_{b\in A'\sm\{a_j\}}\f{x_j-b}{a_j-b}
\\=&\sum_{a_1,\ldots,a_n\in A'}f(a_1,\ldots,a_n)\prod_{j=1}^n\prod_{b\in A'\sm\{a_j\}}\f{x_j-b}{a_j-b}.
\end{align*}
It follows that $\bar f(x_1,\ldots,x_n)$ is a polynomial over the
ring of $p$-adic integers, since $a\not\eq b\ (\mo\ p)$ for any
$a,\,b\in A'$ with $a\not=b$, and $f(a_1,\ldots,a_n)$ are $p$-adic
integers for all $a_1,\ldots,a_n\in A'$ (by the definition of $f$).
As
$$[x_1^{|A|}\cdots x_n^{|A|}]\bar f(x_1,\ldots,x_n)=-c\not\eq0\
(\mo\ p),$$  working in the ring of $p$-adic integers we deduce from
the above that there are $a_1,\ldots,a_n\in A'$ such that
$$f(a_1,\ldots,a_n)\not\eq0\ (\mo\ p).$$

 Note that
$f(0,\ldots,0)=(-1)^{\mathsf{d}^*(G)+p^k-1}-cP(0)^n=0$. So
$I=\{i\in[1,n]:\ a_i\not=0\}$ is nonempty.
 For $i\in I$, we must have $a_i\in A$,  and hence $P(a_i)=0$. It follows that
  $$ \bi{\sum_{i\in[1,n]\sm
I}P(0)-nP(0)-1}{p^k-1}\prod_{t=1}^r\bi{\sum_{s\in
I}a_sc_{st}-1}{p^{k_t}-1}=f(a_1,\ldots,a_n)\not\eq0\ (\mo\ p).$$
With the help of Lemma 2.2, we obtain $-|I|P(0)\eq0\ (\mo\ p^k)$
and\ $\sum_{s\in I}a_sc_{st}\eq0\ (\mo\ p^{k_t})$ for all
$t=1,\ldots,r$.  Therefore $\{c_s\}_{s\in I}$ is an $A$-weighted
zero-sum subsequence of $\{c_i\}_{i=1}^n$ for which $|I|\eq0\ (\mo\
p^k)$ since  $P(0)\not\equiv 0\ (\mo\ p)$.

To see the final part of the theorem, take $k=k_r$ and observe  that
if $(p^k-1)|A|\ge \mathsf{d}^*(G)=\mathsf{D}(G)-1$, then
$$n=p^k-1 +\l\lceil\frac{\mathsf{d}^*(G)+1}{|A|}\r\rceil \le p^k-1 +
\frac{\mathsf{d}^*(G)+|A|}{|A|}\le 2p^k-1,$$
and hence $|I|=p^k$. We are done. \qedsymbol\end{proof}

\section{Terminology and Notation}\label{sec-notation}

In this section, we introduce some more notation to be used in the
remaining part of the paper.

Let $G$ be an abelian group. Then $\Fc(G)$ denotes all  finite, {\it unordered }
sequences (i.e., multi-sets) of $G$ written {\it multiplicatively}.
 We refer to the elements of $\Fc(G)$ as sequences. To lighten the
notation in parts of the paper, we have previously always written sequences with an implicit
order in the format $\{g_i\}_{i=1}^l$, where $g_i\in
G$. However, some of the remaining arguments in
 the paper become more cumbersome to describe without more flexible
notation, so we henceforth use the
{\it multiplicative}
 notation popular among algebraists working in the area (see \cite{Alfred-book} \cite{GG}). In
 particular, a sequence $S\in \Fc(G)$ will be written in the form
 $$S=\prod_{i=1}^{l}g_i=\prod_{g\in G}g^{\vp_g(S)},$$ where $g_i\in G$ are
 the terms in the sequence and  $\vp_g(S)\in \N=\{0,1,2,\ldots\}$ denotes
 the multiplicity of the element $g$ in $S$. Note that the $p$-adic valuation of an integer $x$
 is just the multiplicity of $p$ in the prime factorization of $x=p_1\cdot\ldots\cdot p_l$,
 which is indeed where the notation originates.
 Then $|S|=l$ is the length of
 the sequence, $S'|S$ denotes that $S'$ is a subsequence of $S$ and, in
 such case, ${S'}^{-1}S$ denotes the subsequence of $S$ obtained by
 removing the terms of $S'$ from $S$. The support of $S$, denoted
 $\supp(S)$, consists of all $g\in G$ which occur in $S$, i.e., all $g\in
 G$ with $\vp_g(S)\geq 1$. Of course, if $S,\,T\in \Fc(G)$ are two
 sequences, then $ST\in \Fc(G)$ denotes the sequence obtained by
 concatenating $S$ and $T$. For a homomorphism $\varphi:G\rightarrow G'$,
 we use $\varphi(S)$ to denote the sequence in $G'$ obtained by applying
 $\varphi$ to each term of $S$. Finally, $\sigma(S)=\Sum{i=1}{l}g_i$
 denotes the sum of the terms of the sequence $S$.

 Let $X,\,Y\subseteq G$. Then their sumset is the set $$X+Y=\{x+y\mid x\in
 X,\,y\in Y\}$$ and $-X=\{-x\mid x\in X\}$ denotes the set of inverses of
 $X$. If $A\subseteq \Z$ and $g\in G$, then $$A\cdot g=\{ag\mid a\in
 A\}.$$ We say that $g\in G$ is an $A$-weighted $n$-term subsequence sum
 of $S\in \Fc(G)$, or simply an $A$-weighted $n$-sum of $S$,  if there is
 an $n$-term subsequence $g_1\cdot\ldots\cdot g_n$ of $S$ and $a_i\in A$
 such that $g=\Sum{i=1}{n}a_ig_i$. If we only say $g$ is an $A$-weighted
 subsequence sum of $S\in \Fc(G)$, then we mean it is a $A$-weighted
 $n$-sum of $S$ for some $n\geq 1$. When we say that
 $S=g_1\cdot\ldots\cdot g_n\in \Fc(G)$ has $g$ as an $A$-weighted sum,
this means there are $a_i\in A$ such that $g=\Sum{i=1}{n}a_ig_i$. A
 sequence having the element $0$ as an $A$-weighted sum will simply be
called  an $A$-weighted zero-sum sequence.

 \section{Plus-Minus Weighted Zero-Sums: Generic Bounds and Results for Small
$|G|$}\label{sec-small-values}

In this section,  we focus on $A$-weighted subsequence sums when $A=\{\pm
1\}$ and use the  multiplicative notation for sequences described in
Section \ref{sec-notation}. We begin with an important observation.
Let $G$ be an  abelian group, let $S\in \Fc(G)$ be a sequence, and let
$S'$ be a  subsequence of $S$, say $S'=g_1\cdot\ldots\cdot g_n$ with
$g_i\in G$.  Then $\Sum{i=1}{n}A\cdot g_i$ is the set of all  $A$-weighted
 sums of  $S'$. However, when $A=\{\pm 1\}$, then
$A\cdot  g_i=A\cdot (-g_i)$, and thus the $\{\pm 1\}$-weighted
$n$-term subsequence  sums of $S$ correspond precisely with those of
the sequence  $x^{-1}S(-x)$, for $x\in \supp(S)$ and every $n$. In
other words, we can replace  any term of the sequence $S$ with its
additive inverse without changing which elements of $G$ are
$A$-weighted $n$-term subsequence sums.

When $G$ is an  elementary abelian $2$-group, then $x=-x$ for all $x\in
G$. Consequently,  studying $\{\pm 1\}$-weighted subsequence sums in this
case is no different  than studying ordinary subsequence sums. In
particular (see \cite{Ols0} \cite{Rei}, though the particular cases here are
easy to see),  \be\label{value-22}\mathsf D_{\{\pm 1\}}(\Z_2^2)
=\mathsf  D(\Z_2^2)=3 \und \mathsf s_{\{\pm 1\}}(\Z_2^2)=\mathsf s(\Z_2^2)=5.\ee

The following theorem---and the idea behind its proof---will be one
of the main tools used for proving the results in this section and
the next.

\begin{Theorem}\label{lemma-getzerosums} Let  $G$ be a finite and nontrivial
abelian group and let $S\in \Fc(G)$ be a sequence.

\begin{itemize}

\item[(i)] If $|S|\geq \log_2 |G|+1$ and  $G$ is not an elementary
$2$-group, then $S$ contains a proper,
nontrivial $\{\pm 1\}$-weighted zero-sum subsequence.
    \item[(ii)] If $|S|\geq \log_2 |G|+2$  and $G$ is not an elementary
$2$-group of even rank,
    then $S$ contains a proper, nontrivial  $\{\pm 1\}$-weighted zero-sum
subsequence of even length.
    \item[(iii)] If $|S|>\log_2 |G|$, then $S$ contains a
nontrivial $\{\pm 1\}$-weighted zero-sum  subsequence, and if $|S|>\log_2
|G|+1$, then such a subsequence may be found with even length.
\end{itemize}
\end{Theorem}

\begin{proof}
We begin with the proof of part (i).
Let $S=g_0\cdot g_1\cdot \ldots\cdot g_{l}$, where $g_i\in G$, and set $S'=g_0^{-1}S$.
Note \be\label{stuffty}l=|S'|=|S|-1\geq \log_2 |G|\ee by hypothesis.
There are $2^{l}$ possible subsets  $I\subseteq [1,l]$, each of which corresponds to the sequence
$$S'_I:=\prod_{i\in I}g_i\in \Fc(G)$$ obtained by selecting the terms of $S'$ indexed by the elements of $I$
(including the empty selection $I=\emptyset$, corresponding to the trivial/empty sequence,
which by definition has sum zero).

Suppose there are distinct subsets $I,\,J\subseteq [1,l]$ with
\be\sigma(S'_I)=\label{assump}\Summ{i\in I}g_i=\Summ{j\in
J}g_i=\sigma(S'_J).\ee Since $I\setminus J=I\setminus (I\cap J)$ and
$J\setminus I=J\setminus (I\cap J)$, we can remove the commonly
indexed terms between $S'_I$ and $S'_J$ to find
\be\label{stealth}\sigma(S'_{I\setminus J})=\Summ{i\in I\setminus
J}g_i=\Summ{j\in J\setminus I}g_j=\sigma(S'_{J\setminus I}).\ee
Note, since $I\neq J$, the sets $I\setminus J$ and $J\setminus I$
cannot both be empty, while $I\setminus J$ and $J\setminus I$ are
clearly disjoint. Hence
$$S'_{(I\setminus J)\cup (J\setminus I)}=S'_{I\setminus J}\cdot S'_{J\setminus I}
=\prod_{i\in I\setminus J}g_i\cdot\prod_{j\in J\setminus I}g_j$$
is a nontrivial subsequence of $S'$, which, in view of \eqref{stealth},
has $$0=\Summ{i\in I\setminus J}1\cdot g_i+\Summ{j\in J\setminus I}(-1)\cdot g_j$$
as a $\{\pm 1\}$-weighted sum. Moreover,
since $S'_{(I\setminus J)\cup (J\setminus I)}$ is a subsequence of $S'$ with $S'$ being a proper subsequence of $S$,
 it follows that $S'_{(I\setminus J)\cup (J\setminus I)}$ is a proper $\{\pm\}$-weighted zero-sum subsequence of $S$,
 yielding (i).
 So we may instead assume there do not exist distinct subsets $I,\,J\subseteq [1,l]$ satisfying \eqref{assump}, that is,
 there are no such subsets with $\sigma(S'_I)=\sigma(S'_J)$.

Now \eqref{stuffty} implies that there are $2^{l}\geq |G|$ subsets
$I\subseteq [1,l]$. If $2^l>|G|$, then the pigeonhole principle guarantees
the existence of distinct subsets satisfying \eqref{assump}, contrary to
assumption. Therefore we can assume $2^{l}= |G|$, which is only possible
when equality holds in \eqref{stuffty}: $$|S|= \log_2 |G|+1\in \Z.$$
Moreover, each of the $2^l=|G|$ subsequences $S'_I$, where $I\subseteq
[1,l]$, must have a distinct sum from $G$, else the argument from the
previous paragraph again  completes the proof. In consequence,
every element of $G\setminus \{0\}$ is  representable as a subsequence sum of $S'$ with
$0$ represented by the  trivial sequence. In particular, it follows that
there exist subsequences  $T_1$ and $T_2$ of $S'$ with $\sigma(T_1)=g_0$
and $\sigma(T_2)=-g_0$,  where (recall) $g_0$ is the term from $S$ that we
removed to obtain $S'$.  Consequently, if $T_1$ is a proper subsequence of
$S'$, then $g_0T_1$ is a  proper $\{\pm 1\}$-weighted zero-sum subsequence
of $S$, while if $T_2$ is  a proper subsequence of $S'$, then $g_0T_2$ is
a proper $\{\pm 1\}$-weighted  zero-sum subsequence of $S$. In either
case, the proof of (i) is  complete, so we must have $S'=T_2=T_1$, in
which case  \be\label{one}\sigma(S')=-g_0=\sigma(T_2)=\sigma(T_1)=g_0.\ee

Since  every element of $G$  occurs as the sum of one of the $|G|$
subsequences $S'_I$, where $I\subseteq [1,l]$, and since
$I=\emptyset$ corresponds to the subsequence with sum $0$, we
conclude that $\la \supp(S')\ra=G$. Consequently, since $G$ is not
an elementary $2$-group, it follows that there must be some $y\in
\supp(S')$ with \be\label{two}2y\neq 0.\ee Now, recall that
replacing a term from a sequence with its additive inverse does not
affect any of the $\{\pm 1\}$-weighted subsequence sums (as
explained at the beginning of the section). Thus, it suffices to
prove (i) for the sequence $S_0:=y^{-1}S(-y)$ obtained by replacing
$y$ by $-y$ in $S$. Note that $$\sigma(S_0)=\sigma(S)-2y\und
\sigma(S'_0)=\sigma(S')-2y,$$ where $S_0':=y^{-1}S'(-y)$. Therefore,
using \eqref{one}, we derive that
\be\label{three}\sigma(S'_0)=\sigma(S')-2y=g_0-2y.\ee However,
applying all of the above arguments using $S_0=y^{-1}S(-y)$ and
$S_0'=y^{-1}S'(-y)$, we will complete the proof unless \eqref{one}
holds for $S'_0$ as well:  $$\sigma(S'_0)=g_0.$$ Combining this
equality with \eqref{three}, we find that $2y=0$, which contradicts
\eqref{two}, completing the proof of (i).

\medskip

We continue with the proof of part (ii), which is just a variation
on that of (i). As before, let $S=g_0\cdot g_1\cdot \ldots\cdot
g_{l}$, where $g_i\in G$, and set $S'=g_0^{-1}S$. Note that
\be\label{stuffty-p} l=|S'|=|S|-1\geq \log_2 |G|+1\ee by hypothesis.
By a well-known combinatorial identity (which can be proven using a
simple inductive argument and the correspondence between a subset
and its compliment), there are
$$2^{l-1}=\Sum{i=0}{\lfloor
l/2\rfloor}\binom{l}{2i}=\Sum{i=0}{\lfloor
(l-1)/2\rfloor}\binom{l}{2i+1}$$ possible subsets  $I\subseteq
[1,l]$ of odd cardinality, each of which corresponds to the odd
length sequence
$$S'_I:=\prod_{i\in I}g_i\in \Fc(G)$$ obtained by selecting the terms of $S'$ indexed by the elements of $I$.

Suppose there are distinct subsets $I,\,J\subseteq [1,l]$ of odd
cardinality with \be\sigma(S'_I)=\label{assump-p}\Summ{i\in
I}g_i=\Summ{j\in J}g_i=\sigma(S'_J).\ee Since $I\setminus
J=I\setminus (I\cap J)$ and $J\setminus I=J\setminus (I\cap J)$, we
can remove the commonly indexed terms between $S'_I$ and $S'_J$ to
find   \be\label{stealth-p}\sigma(S'_{I\setminus J})=\Summ{i\in
I\setminus J}g_i=\Summ{j\in J\setminus I}g_j=\sigma(S'_{J\setminus
I}).\ee Note, since $I\neq J$, the sets $I\setminus J$ and
$J\setminus I$ cannot both be empty, while $I\setminus J$ and
$J\setminus I$ are clearly disjoint; furthermore, $|I\setminus
J|+|J\setminus I|=|I|+|J|-2|I\cap J|$ is an even number in view of
$|I|\equiv |J|\ (\mo\ 2)$. Hence
$$S'_{(I\setminus J)\cup (J\setminus I)}=S'_{I\setminus J}\cdot S'_{J\setminus I}
=\prod_{i\in I\setminus J}g_i\cdot\prod_{j\in J\setminus I}g_j$$ is a nontrivial subsequence of $S'$ with even length,
which, in view of \eqref{stealth-p}, has $$0=\Summ{i\in I\setminus J}1\cdot g_i+\Summ{j\in J\setminus I}(-1)\cdot g_j$$
as a $\{\pm 1\}$-weighted sum. Moreover, since $S'_{(I\setminus J)\cup (J\setminus I)}$ is a subsequence of $S'$
with $S'$ being a proper subsequence of $S$, it follows that $S'_{(I\setminus J)\cup (J\setminus I)}$ is
a proper $\{\pm 1\}$-weighted zero-sum subsequence of $S$, yielding (ii).
So we may instead assume there do not exist distinct subsets $I,\,J\subseteq [1,l]$ of odd cardinality satisfying
\eqref{assump-p}, that is, there are no such subsets with $\sigma(S'_I)=\sigma(S'_J)$.

Now \eqref{stuffty-p} implies that there are $2^{l-1}\geq |G|$ subsets
$I\subseteq [1,l]$ of odd cardinality. If $2^{l-1}>|G|$, then the pigeonhole
principle guarantees the existence of distinct subsets of odd cardinality
satisfying \eqref{assump-p}, contrary to assumption. Therefore we can
assume $2^{l-1}= |G|$, which is only possible when equality holds in
\eqref{stuffty-p}: \be\label{seeal}|S'|= \log_2 |G|+1\in \Z.\ee Moreover,
each of the $2^{l-1}=|G|$ odd length subsequences $S'_I$ must have a
 distinct sum from $G$, else the argument from the previous paragraph
again  completes the proof. In consequence, every element of $G$ is
representable as an odd length subsequence sum of $S'$. In particular,
it follows that there exist odd length subsequences $T_1$ and $T_2$ of $S'$
with $\sigma(T_1)=g_0$ and $\sigma(T_2)=-g_0$, where (recall) $g_0$ is the term from $S$ that we removed to obtain $S'$.
Consequently, if $T_1$ is a proper subsequence of $S'$, then $g_0T_1$ is a proper $\{\pm 1\}$-weighted
zero-sum subsequence of $S$ of even length
(since the length of $T_1$ is odd), while if $T_2$ is a proper subsequence of $S'$, then $g_0T_2$
is a proper $\{\pm 1\}$-weighted zero-sum
subsequence of $S$ of even length (since the length of $T_2$ is odd). In either case, the proof of (ii) is complete,
so we must have $S'=T_2=T_1$
with $|S'|=|T_1|=|T_2|$ odd, whence  \be\label{one-p}\sigma(S')=-g_0=\sigma(T_2)=\sigma(T_1)=g_0.\ee

If $G$ is an elementary $2$-group, then $\log_2 |G|$ is the rank of $G$, which is assumed odd by hypothesis.
But in this case, \eqref{seeal}
implies that  $|S'|=\log_2 |G|+1$ is an even number, contrary to what we have just seen above.
Therefore we may assume $G$ is not
an elementary $2$-group.

Since  every element of $G$ occurs as the sum of one of the $|G|$ odd length subsequences $S'_I$ of $S'$,
we conclude that $\la \supp(S')\ra=G$.
Consequently, since $G$ is not an elementary $2$-group, it follows that there must be some $y\in \supp(S')$
 with \be\label{two-p}2y\neq 0.\ee
The remainder of the proof now concludes identical to that of part (i).
%

\medskip

The proof of part (iii) is a routine simplification of the proofs of parts (i) and (ii).\qedsymbol
\end{proof}

\bigskip

Next, we give the proof of Theorem \ref{cor-basic-bounds-dav}, which
is a simple corollary of Theorem \ref{lemma-getzerosums}.

\medskip

\begin{proof}[Proof of Theorem \ref {cor-basic-bounds-dav}] We begin with
the first set of inequalities. The upper bound follows from Theorem
\ref{lemma-getzerosums}(iii). We turn to the lower bound.

Let $e_1,\ldots,e_r$ be a basis for $G$ with $G=\langle
e_1\rangle\oplus \langle e_2\rangle\oplus\cdots\oplus \langle
e_r\rangle$ and  $\langle e_i\rangle \cong \Z_{n_i}$ for
$i=1,\ldots,r$. For $i\in [1,r]$, define  $$S_i=(2^0e_i)\cdot
(2^1e_i)\cdot (2^2e_i)\cdot\ldots\cdot (2^{\lfloor \log_2
n_i\rfloor-1}e_i)\in \Fc(\langle e_i\rangle)$$ and then set
$$S=S_1S_2\cdot\ldots\cdot S_r\in\Fc(G).$$ Note
$|S|=\Sum{i=1}{r}\lfloor\log_2 n_i\rfloor$. Thus it suffices to show
that $S$ contains no nontrivial $\{\pm 1\}$-weighted zero-sum
subsequence. Moreover, since the $e_1,\ldots,e_r$ form a basis of
$G$, it in fact suffices to show that each $S_j$, for
$j=1,\ldots,r$, contains no nontrivial $\{\pm 1\}$-weighted zero-sum
subsequence.

Let $j\in [1,r]$ and consider an arbitrary $\{\pm 1\}$-weighted
subsequence sum of $S_j$, say $$\Sum{i=0}{\lfloor \log_2
n_j\rfloor-1}\varepsilon_i 2^ie_j\quad\ \mbox{ with } \{0\}\not=\{\varepsilon_i:\,i\ge0\}\se\{0,\pm1\}.$$ We will show
that $\Sum{i=0}{\lfloor \log_2 n_j\rfloor-1}\varepsilon_i 2^ie_j\neq
0$. Let $t\in [0,\lfloor \log_2 n_j\rfloor-1]$ be the maximal index
such that $\varepsilon_t\neq 0$ and w.l.o.g. assume
$\varepsilon_t=1$ (by multiplying all terms by $-1$ if necessary).
Then
$$0<1= 2^t-\Sum{i=0}{t-1}2^i\leq  2^t+\Sum{i=0}{t-1}\varepsilon_i 2^i= \Sum{i=0}{\lfloor \log_2 n_j\rfloor-1}
\varepsilon_i 2^i \leq \Sum{i=0}{\lfloor \log_2 n_j\rfloor-1}2^i
\leq n_j-1,$$ which shows that $\Sum{i=0}{\lfloor \log_2
n_j\rfloor-1}\varepsilon_i 2^ie_j\neq 0$ (since $\ord(e_j)=n_j$).
Consequently, $S_j$ contains no $\{\pm 1\}$-weighted zero-sum
subsequence, for each $j\in [1,r]$, showing that
$$\mathsf D_{\{\pm 1\}}(G)\geq |S|+1=\Sum{i=1}{r} \lfloor\log_2 n_i\rfloor+1.$$

To show the second set of inequalities, let $S\in \Fc(G)$ be a sequence of length $|S|=\mathsf D_{\{\pm
 1\}}(G)-1$ containing no $\{\pm 1\}$-weighted zero-sum subsequence. It is
 then clear that the sequence $0^{n_r-1}S$ contains no $\{\pm
 1\}$-weighted zero-sum subsequence of length $n_r$, showing the first
inequality, while the second inequality follows by the first part.\qedsymbol
\end{proof}

\medskip

To show that our method can also be used to precisely determine $\mathsf s_{\{\pm 1\}}(G)$ in certain cases,
we will compute the values of $\mathsf s_{\{\pm 1\}}(\Z_4^2)$ and $\mathsf s_{\{\pm 1\}}(\Z_8^2)$.
These will then be used to give a simple bound for
$\mathsf s_{\{\pm 1\}}(\Z_n^2)$ complementing the result of \cite{ABPP}. We begin with the following lemma.

\begin{Lemma}\label{lem-hypgraph-case8} Let $G=\Z_8^2$ and let $S\in \Fc(G)$ be a sequence with $|S|=10$.
Then $S$ contains a $\{\pm 1\}$-weighted zero-sum subsequence $T$ of length $|T|\in \{2,4,8\}$.
\end{Lemma}

\begin{proof}
We adapt the proof of Theorem \ref{lemma-getzerosums}. Let
$S=g_1\cdot g_2\cdot \ldots\cdot g_{10}$, where $g_i\in G$. For
$j\in [0,10]$, let $\mathcal{I}_j$ be the set of all subsets
$I\subseteq [1,10]$ having cardinality $j$. We consider
$X:=\mathcal{I}_4\cup \mathcal{I}_2$. Recall that we associate each
$I\subseteq [0,10]$ with the indexed subsequence $S_I:=\prod_{i\in
I}g_i$ of $S$.  We now analyze the possible intersection
cardinalities between sets $I,\,J\in X$ with
\be\label{eq-holds}\sigma(S_I)=\Summ{i\in I}g_i=\Summ{j\in
J}g_j=\sigma(S_J).\ee

Let $I,\,J\in X$ be distinct indexing subsets such that \eqref{eq-holds} holds. Thus, by removing terms contained in both
$S_I$ and $S_J$, we obtain $$\sigma(S_{I\setminus J})
=\Summ{i\in I\setminus J}g_i=\Summ{j\in J\setminus I}g_j=\sigma(S_{J\setminus I}).$$
 Note, since $I\neq J$, the sets $I\setminus J$ and $J\setminus I$ cannot both be empty,
 while clearly  $I\setminus J$ and $J\setminus I$ are disjoint.
 Hence $$S_{(I\setminus J)\cup (J\setminus I)}=S_{I\setminus J}\cdot S_{J\setminus I}
 =\prod_{i\in I\setminus J}g_i\cdot\prod_{j\in J\setminus I}g_j\in \Fc(G)$$ is a nontrivial subsequence of $S$ having
 $$0=\Summ{i\in I\setminus T}1\cdot g_i+\Summ{j\in J\setminus T}(-1)\cdot g_j$$ as a $\{\pm 1\}$-weighted zero-sum.
 Assuming by contradiction that
 $S$ contains no $\{\pm 1\}$-weighted zero-sum subsequence $T$ of length $|T|\in \{2,4,8\}$,
 we conclude that \be|S_{(I\setminus J)\cup (J\setminus I)}|=|I\setminus J|+|J\setminus I|=|I|+|J|-2|I\cap J|\notin
 \{2,4,8\}.\nn\ee
 Using the above restriction, for distinct indexing sets $I,\,J\in X$ satisfying
 \eqref{eq-holds}
 we see that $|I|=|J|=2$ is impossible, and
 \be\label{card-restric-transfred} |I\cap J|
=\left\{
   \begin{array}{ll}
     0, & \hbox{if } |I|=2 \und |J|=4,\\
     1, & \hbox{if } |I|=4\und |J|=4.,
   \end{array}
 \right.\ee

Using \eqref{card-restric-transfred}, we proceed to estimate the maximal number of subsets $I\in X$ that
can simultaneously
have all their corresponding subsequences $S_I$ being of equal sum. Observe these are just very
particular $L$-intersecting set system problems
over $|S|=10$ vertices. To this end, let $I_1,\ldots,I_n\in X$ be distinct indexing subsets
with $\sigma(S_{I_j})=\sigma(S_{I_k})$ for all $j$ and $k$.
We proceed with some useful comments regarding the $I_j$ under this assumption of equal sums.

\begin{itemize}

\item In view of \eqref{card-restric-transfred}, there can be at most one  $I_j$ with $|I_j|=2$.

\item If (say) $|I_1|=|I_2|=|I_3|=4$ with $|I_1\cap I_2\cap I_3|=1$, then \eqref{card-restric-transfred}
 implies $|I_1\cup I_2\cup I_3|=10=|S|$. Thus,
since any further $I_j$ with $|I_j|=2$ must be disjoint from any other $I_k$
(in view of  \eqref{card-restric-transfred}),
 and since $|S|=10=|I_1\cup I_2\cup I_3|$, we see in this case that no $I_j$ has $|I_j|=2$.
 Therefore, if there is a further $I_j$ with $j\geq 4$,
  then it must have cardinality $|I_j|=4$, in which case \eqref{card-restric-transfred} shows that $I_j$ can contain
  at most one element from each
  set $I_1$, $I_2$ and $I_3$. However, since $|I_1\cup I_2\cup I_3|=10=|S|$, there are no further elements to be found,
  whence $|I_j|\leq 3$, a contradiction.

In summary, if three sets $I_j$ of size $4$ intersect in a common point, then $n=3$ and there are no other
indexing sets $I_j$ besides these three.

\item If (say) $|I_1|=|I_2|=|I_3|=4$ with $|I_1\cap I_2\cap I_3|\neq 1$, then \eqref{card-restric-transfred}
ensures that these sets lie as depicted in the following diagram,
where each line below represents one of the sets $I_j$,
where $j\in [1,3]$, with the points contained in the line corresponding to the elements of $I_j$.
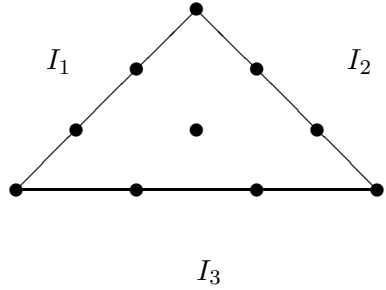
\begin{figure}[h]\label{caseone}
\setlength{\unitlength}{8mm}
\begin{center}
\begin{picture}(6,5)
\put(1,1){\circle*{0.2}}
\put(2,2){\circle*{0.2}}
\put(3,3){\circle*{0.2}}
\put(4,4){\circle*{0.2}}
\put(4,2){\circle*{0.2}}
\put(3,1){\circle*{0.2}}
\put(5,1){\circle*{0.2}}
\put(7,1){\circle*{0.2}}
\put(6,2){\circle*{0.2}}
\put(5,3){\circle*{0.2}}

\put(4,-0.5){$I_3$}
\put(6.5,3){$I_2$}
\put(1.5,3){$I_1$}

\put(1,1){\line(1,1){3}}
\put(1,1){\line(1,0){6}}
\put(7,1){\line(-1,1){3}}
\end{picture}
\end{center}
\vspace{3mm} \caption{Configuration for $3$ intersecting $4$-sets with no common intersection}
\end{figure}

Since \eqref{card-restric-transfred} ensures that any $I_j$ with $|I_j|=2$ must be disjoint from all other $I_i$,
we see there can be no such $I_j$ in this case. Using \eqref{card-restric-transfred} and the previous comment,
it is now easily verified that there can be at most two additional $I_j$ with $|I_j|=4$ besides $I_1$, $I_2$ and $I_3$
(as each of the new points from these additional $I_j$, for $j\geq 4$, must avoid points already covered by two edges,
 such as the three corners of the triangle depicted above).

In summary, if there are three $I_j$ of size $4$ that do not intersect in a common point, then no $I_j$
has size $|I_j|=2$ and $n\leq 5$.

\item In view of the previous remarks, we see that if some $I_j$ has $|I_j|=2$, then $n\leq 3$ and all other $I_k$
with $k\neq j$
have $|I_k|=4$: the first remark ensures that all other $I_k$ have $|I_k|=4$, while the second and third remark combine
 to imply there are at most two $I_i$ with $|I_i|=4$.

\item Combining the last three comments, we see that if no $I_j$ has $|I_j|=2$, that is, $|I_j|=4$ for all $j$,
then $n\leq 5$.
\end{itemize}

There are $\binom{10}{4}=210$ subsets $I\in \mathcal I_4$ and $\binom{10}{2}=45$ subsets $I\in \mathcal I_2$. Clearly,
many of their corresponding sequences $S_I$ must have common sum as there are only $|G|=64$ sums to choose from.
It is clear, from the final two comments above, that in order to minimize the number of sums spanned by
all $I\in X=\mathcal I_4\cup \mathcal I_2$,
we must pair each $J\in \mathcal I_2$ with two $I,\,I'\in \mathcal I_4$ (to form a grouping corresponding to some
distinct sum from $G$)
and then take all remaining (unpaired) $I\in \mathcal I_4$  and put them into groupings of $5$
(with one leftover remainder group possible,
and each of these groupings corresponding to some distinct sum from $G$). In other words, there are at least
$45+\frac15(210-2\cdot 45)=69$ distinct
 sums covered by the the sets $I\in X$. However, since there are only $|G|=64<69$ sums available,
 this is a contradiction, completing the
 proof (essentially, one of the intersection conditions given by \eqref{card-restric-transfred} must actually fail,
 which then gives rise
 to the weighted
  zero-sum subsequence of one of the desired lengths).\qedsymbol
\end{proof}

As promised, we now compute the values of $\mathsf s_{\{\pm 1\}}(\Z_4^2)$ and $\mathsf s_{\{\pm 1\}}(\Z_8^2)$
and use them to give a simple upper bound for $\mathsf s_{\{\pm 1\}}(\Z_{n}^2)$. The method below could also be iterated
to obtain progressively better bounds for larger $u=\vp_2(n)$. However, in view of the results of the next section,
we do not expand upon this.
\begin{Theorem} \label{thm-constant-2-epsilon}Let $n\in \Z^+$. Then $\mathsf s_{\{\pm 1\}}(\Z_{n}^2)\leq 2n+1$.
Indeed, letting $u=\vp_2(n)$ denote the maximum power of $2$ dividing $n$, we have the following bounds:
\begin{itemize}
\item[(i)] $\mathsf s_{\{\pm 1\}}(\Z_2^2)=5$, $\mathsf s_{\{\pm 1\}}(\Z_4^2)=8$ and $\mathsf s_{\{\pm 1\}}(\Z_8^2)=14$.
\item[(ii)] If $u\leq 1$, then $\mathsf s_{\{\pm 1\}}(\Z^2_{n})\leq 2n+1$.
\item[(iii)] If $u=2$, then $\mathsf s_{\{\pm 1\}}(\Z^2_{n})\leq 2n$.
\item[(iv)] If $u\geq 3$, then $\mathsf s_{\{\pm 1\}}(\Z_{n}^2)\leq  \frac{15}{8}n+1$.
\end{itemize}
\end{Theorem}

\begin{proof}[Proof of Theorem \ref{thm-constant-2-epsilon}]

As mentioned in the introduction, if $n$ is odd, then
 \be\label{stardust}\mathsf s_{\{\pm 1\}}(\Z_n^2)=2n-1\leq
2n+1.\ee Thus we restrict our attention to the case $2\mid m$.

We proceed to prove part (i), which contains the crucial basic cases used in the  inductive approach.
The case $\Z_2^2$ is
covered by \eqref{value-22}, so we begin with $\Z_4^2$.

By Theorem \ref{cor-basic-bounds-dav}, we have $\mathsf s_{\{\pm
 1\}}(\Z_4^2)\geq 4+2\cdot 2=8$.  It remains to show that $\mathsf
 s_{\{\pm 1\}}(\Z_4^2)\leq 8$. Note that $\exp(\Z_4^2)=4$ and
 $\log_2 |\Z_4^2| =4$. Let $S\in \Fc(\Z_4^2)$ be a sequence with
 length $|S|=8$. Applying Theorem \ref{lemma-getzerosums}(ii) to a
 subsequence of $S$ of length $6$, we obtain a $\{\pm 1\}$-weighted
 zero-sum subsequence $T$ of $S$ with length $|T|\in \{2,4\}$. We may
 assume $|T|=2$, else the desired length weighted zero-sum subsequence is found. But now, applying Theorem
 \ref{lemma-getzerosums}(ii) to $T^{-1}S$, we likewise obtain another
 $\{\pm 1\}$-weighted zero-sum subsequence $T'$ of $T^{-1}S$ with length
 $|T'|=2$, and then $TT'$ is a $\{\pm 1\}$-weighted zero-sum subsequence
of $S$ with length $4$. This shows that $\mathsf s_{\{\pm
 1\}}(\Z_4^2)\leq 4+2\cdot 2=8$

We continue with the case $\Z_8^2$.
 By Theorem \ref{cor-basic-bounds-dav}, we have $\mathsf s_{\{\pm
 1\}}(\Z_8^2)\geq 8+2\cdot 3=14$.  It remains to show that $\mathsf
 s_{\{\pm 1\}}(\Z_8^2)\leq 14$. Note that $\exp(\Z_8^2)=8$  and
 $\log_2 |\Z_8^2| =6$. Let $S\in \Fc(\Z_8^2)$ be a sequence with
 length $|S|=14$ and assume by contradiction that $S$ contains no $\{\pm
 1\}$-weighted zero-sum subsequence of length $8$.

\medskip

{\it Case 1: } There exists a $\{\pm 1\}$-weighted zero-sum subsequence
$S_0$ of $S$ with length $|S_0|=2$.
 Let $R$ be a maximal length subsequence of $S$ such that $|R|$ is even
 and, for every $n\in [0,|R|]\cap 2\Z$, $R$ contains a $\{\pm
 1\}$-weighted zero-sum subsequence $T_n$ of length $|T_n|=n$. In view of
 the case hypothesis, $R$ exists with $|R|\geq 2$.
 Note that $|R|\leq 6$, else $R$, and hence also $S$, contains a weighted
 zero-sum subsequence of length $8$, as desired. Thus
 $|R^{-1}S|=|S|-|R|\geq 14-6=8$. Applying Theorem
 \ref{lemma-getzerosums}(ii) to a subsequence of $R^{-1}S$ of length $8$,
 we obtain a $\{\pm 1\}$-weighted zero-sum subsequence $T$ of $R^{-1}S$
with length $|T|\in \{2,4,6\}$.

 Suppose $|T|\leq |R|+2$.  Then define $R'=RT$. For every $n\in
 [2,|R|]\cap 2\Z$, we see that $R$, and hence also $R'$, contains the
 $\{\pm 1\}$-weighted zero-sum subsequence $T_n$ of length $|T_n|=n$. On
 the other hand, since $|T|\leq |R|+2$ and $|R|$ is even, it follows that
 every $m\in [|R|+2,|R|+|T|]\cap 2\Z$ can be written in the form
 $$m=|T|+n\;\mbox{ with }\;n\in [|R|-|T|+2,|R|]\und n\geq 0.$$ Thus, since $|T|$ is even, it follows that
 the
 subsequence $T_nT$ of $RT$ is a weighted zero-sum subsequence of length $m\in
 [|R|+2,|R|+|T|]\cap 2\Z$, in which case $RT$ contradicts the maximality
of $R$. So we conclude that $|T|\geq |R|+4\geq 6$.

 Hence, since $|T|\in \{2,4,6\}$, we see that $T|R^{-1}S$ is a weighed
zero-sum subsequence of length $|T|=6$.
 But now, since $|R|\geq 2$, it follows, by the defining property of $R$,
 that there exists  $\{\pm 1\}$-weighted zero-sum subsequence $T_2|R$ with
 length $|T_2|=2$, whence $T_2T$ is a $\{\pm 1\}$-weighted zero-sum
 subsequence of $S$ with length $|T_2|+|T|=2+6=8$, as desired. This
completes the case.

\medskip

{\it Case 2: } There does not exist a $\{\pm 1\}$-weighted zero-sum
subsequence $S_0$ of $S$ with length $|S_0|=2$.
 Since we have assumed by contradiction that $S$ contains no weighted
 zero-sum subsequence of length $8$, and in view of the hypothesis of the case, we see
 that applying Lemma \ref{lem-hypgraph-case8} to a subsequence of $S$ of
 length $10$ yields a $\{\pm 1\}$-weighted zero-sum subsequence $T$ of $S$
 with $|T|=4$. Noting that $|T^{-1}S|=10$, we see that a second
 application of Lemma \ref{lem-hypgraph-case8} to $T^{-1}S$ yields another
 $\{\pm 1\}$-weighted zero-sum subsequence $T'|T^{-1}S$ with $|T'|=4$. But
 now $TT'$ is a $\{\pm 1\}$-weighted zero-sum subsequence of $S$ with
 length $|T|+|T'|=4+4=8$, as desired, which completes the proof of part (i).

 \bigskip

Next we prove part (ii). In view of part (i), we have $\mathsf s_{\{\pm
 1\}}(\Z_{2}^2)=5=2n+1$. Thus, in view of \eqref{stardust},  we
 may assume $n=2m$ with $m>1$ odd. Let $\varphi: \Z_{2m}^2\rightarrow m\cdot\Z_{2m}^2$
 denote the multiplication by $m$ homomorphism, which has kernel $\ker
\varphi\cong \Z_m^2$ and image $\varphi(\Z_{2m}^2)\cong \Z_2^2$.
 Note that $$|S|=|\varphi(S)|=2n+1=2(2m-2)+5.$$ Thus, iteratively applying the
 definition of $\mathsf s_{\{\pm 1\}}(\Z_{2}^2)=5$ to the
 sequence $\varphi(S)$, we find $2m-1$ subsequences
 $S_1,\ldots,S_{2m-1}\in \Fc(\Z_{2m}^2)$, each of length $|S_i|=2$, such that
 $S_1S_2\cdot \ldots\cdot S_{2m-1}|S$ and each $S_i$ has a $\{\pm
 1\}$-weighted sum $x_i\in \ker \varphi\cong \Z_m^2$. Observe, by
 swapping the signs on every term of $S_i$, that $-x_i\in  \ker
 \varphi\cong \Z_m^2$ is also a $\{\pm 1\}$-weighted sum of
 $S_i$. Now applying the definition of $\mathsf s_{\{\pm 1\}}(\Z_m^2)=2m-1$
 (see \eqref{stardust}) to the sequence $x_1\cdot \ldots\cdot
 x_{2m-1}$, we find an $m$-term subsequence, say $x_1\cdot\ldots \cdot
 x_m$, having $0$ as a $\{\pm 1\}$-weighted sum. However, since each $\pm
 x_i$ was a  $\{\pm 1\}$-weighted sum of the subsequence $S_i$, we
 conclude that the subsequence $S_1\cdot \ldots\cdot S_m|S$ has $0$ as a
 $\{\pm 1\}$-weighted sum. Since each $S_i$ has length $2$, we see
 $|S_1\cdot \ldots\cdot S_m|=2m=n$. Thus $S_1\cdot \ldots \cdot S_m$ is a
 $\{\pm 1\}$-weighted zero-sum subsequence of $S$ with length $n$, as desired.

 \bigskip

Next, the proof of part (iii). Let $n=4m$ with $m$ odd. By part (i), we know
 $\mathsf s_{\{\pm 1\}}(\Z_4^2)=8=2n$. Thus we may assume $m>1$.
 Let $\varphi: \Z_{4m}^2\rightarrow m\cdot \Z_{4m}^2$ denote the multiplication by $m$
 homomorphism, which has kernel $\ker \varphi\cong \Z^2_m$ and
image $\varphi(\Z_{4m}^2)\cong \Z_4^2$.
 Note that $$|S|=|\varphi(S)|=2n=4(2m-2)+8.$$ Thus, iteratively applying the
 definition of $\mathsf s_{\{\pm 1\}}(\Z_{4}^2)=8$ to the
 sequence $\varphi(S)$, we find $2m-1$ subsequences
 $S_1,\ldots,S_{2m-1}\in \Fc(\Z_{4m}^2)$, each of length $|S_i|=4$, such that
 $S_1S_2\cdot \ldots\cdot S_{2m-1}|S$ and each $S_i$ has a $\{\pm
 1\}$-weighted sum $x_i\in \ker \varphi\cong \Z^2_m$. Then, as in
 part (ii), applying the definition of $\mathsf s_{\{\pm 1\}}(\Z_m^2)=2m-1$
 (see \eqref{stardust}) to the sequence $x_1\cdot \ldots\cdot
 x_{2m-1}$ yields a subsequence (say) $S_1\cdot \ldots\cdot S_m|S$ which
 is a $\{\pm 1\}$-weighted zero-sum subsequence of length
$|S_1\cdot\ldots\cdot S_m|=4m=n$, as desired.

\bigskip

Finally, we conclude with the proof of part (iv). Let $n=8m$ with
$m\in \Z^+$. If $n=8$, then part (i) implies $\mathsf s_{\{\pm
1\}}(\Z^2_8)=14<\frac{15}{8}n+1=16$. We proceed by induction on $n$.
Let $S\in \Fc(\Z_{8m}^2)$ be a sequence  with $|S|\geq
\frac{15}{8}n+1=15m+1$. Let $\varphi:G\rightarrow 8\cdot  G$ be the
multiplication by $8$ map, which has kernel $\ker \varphi\cong
\Z_{8}^2$ and image $\varphi(\Z_{8m}^2)=8 \cdot \Z_{8m}^2\cong
\Z_m^2$. By induction hypothesis or parts (ii) and (iii), we
conclude that $\mathsf s_{\{\pm\}}(\Z_m^2)\leq 2m+1$. Note that
$$|S|=|\varphi(S)|=13m+(2m+1).$$ Thus, iteratively applying the
 definition of $\mathsf s_{\{\pm 1\}}(\Z_{m}^2)\leq 2m+1$ to the
 sequence $\varphi(S)$, we find $14$ subsequences
 $S_1,\ldots,S_{2m-1}\in \Fc(\Z_{8m}^2)$, each of length $|S_i|=m$, such that
 $S_1S_2\cdot \ldots\cdot S_{14}|S$ and each $S_i$ has a $\{\pm
 1\}$-weighted sum $x_i\in \ker \varphi\cong \Z^2_8$. Then, as in
 parts (ii) and (iii), applying the definition of $\mathsf s_{\{\pm 1\}}(\Z_8^2)=14$ (from part (i))
 to the sequence $x_1\cdot \ldots\cdot
 x_{14}$ yields a subsequence (say) $S_1\cdot \ldots\cdot S_8|S$ which
 is a $\{\pm 1\}$-weighted zero-sum subsequence of length
$|S_1\cdot\ldots\cdot S_8|=8m=n$, as desired.\qedsymbol
\end{proof}

\section{Plus-Minus Weighted Zero-Sums: Asymptotic Bounds}

For the proof of Theorem \ref{thm-even-assymptotics}, we will need
to make use of several results from the theory of $L$-intersecting
set systems. The following is now a well-known result from this
area. See \cite{Ray-chaudhur-Wilson} for the original $t$-design
formulation, \cite{Frannkl-Wilson} for a more general mod $p$
formulation, and \cite{ABS} for a yet more general result.

\begin{Theorem}[Uniform Frankl-Ray-Chaudhuri-Wilson Theorem]\label{Ray-Chaudhuri-Wilson-thm} Let $k,\,n\in \Z^+$
be integers,
let $\mathcal F$ be a collection of $k$-element subsets of an $n$-element set, and let $L\subseteq \{0,1,2,\ldots,k-1\}$
be a subset.
Suppose \be\nn |E\cap E'|\in L\quad \mbox{ for all distinct }\; E,\,E'\in \mathcal F.\ee
Then $|\mathcal F|\leq \binom{n}{|L|}$.
\end{Theorem}

We will also need a more recent prime power version of the
Nonuniform Frankl-Ray-Chaudhuri-Wilson Inequality
\cite{prime-power-frankl-wilson}. To state it, we must first
introduce the following definition. We say that a polynomial
$f(x)\in \Z[x]$ {\it separates} the element $\alpha \in \Z$ from the
set $B\subseteq \Z$ with respect to the prime $p$ if $$\vp_p
(f(\alpha))<\min_{b\in B}\;\vp_p(f(b)),$$ where $\vp_p(x)$ denotes
the $p$-adic valuation of a rational number $x$ (and $\vp_p(0)$ is
regarded as $+\infty$).

\begin{Theorem}\label{thm-primepower} Let $p$ be a prime, let $q=p^k$ with $k\geq 1$, and let $K$ and $L$
be disjoint subsets of $\{0,1,\ldots,q-1\}$. Let $\mathcal F$ be a collection of subsets of an $n$-element set.
Suppose \ber\nn|E|&\in& K+q\Z\;\mbox{ for all }\;E\in \mathcal F\und\\ |E\cap E'|&\in& L+q\Z\;
\mbox{ for all distinct }\;E,\,E'\in \mathcal{F}.\nn\eer
Then $|\mathcal F|\leq \binom{n}{D}+\binom{n}{D-1}+\ldots+\binom{n}{0}$, where $D\leq 2^{|L|-1}$
is the maximum over all $\alpha\notin L$
of the minimal degree of a polynomial separating the element $\alpha$ from the set $L+q\Z$   with respect to $p$.
\end{Theorem}

For $m\in\Z^+$, $n\in\N$ and $r\in\Z$ we set
$$\sbinom{n}{r}_m=\Sum{\underset{i\equiv r\,({\rm mod}\ m)}{0\le i\le n}}{}\binom{n}{i}.$$

\begin{Lemma} {\rm (\cite[Remark 1.1]{Sun08})}\label{sbinom-lemma} For any $m,n\in\Z^+$, we have
$$\sbinom{n}{\lfloor \frac{n+1}{2}\rfloor}_m\geq
\frac{2^n}{m},$$
and furthermore
\be\label{goaly}\sbinom{n}{\lfloor\frac{n+1}{2}\rfloor}_m\geq
\sbinom{n}{\lfloor\frac{n+1}{2}\rfloor+1}_m \geq \cdots\geq
\sbinom{n}{\lfloor\frac{n+m}{2}\rfloor}_m.\ee
\end{Lemma}

With the above tools in hand, we can conclude the proof of Theorem
\ref{thm-even-assymptotics}.
 With regards to asymptotic notation, recall that $f(x)=O(g(x))$ (or $f(x)\ll g(x)$) (as $x\to+\infty$)
  means that there exists a constant $C> 0$
 such that $|f(x)|\leq C|g(x)|$ for all sufficiently large values of $x$, while $f(x)\gg g(x)$
 means that there exists a constant $C>0$ such that $|f(x)|\geq C|g(x)|$ for all sufficiently large values of $x$,
 where $f$ and $g$ are functions.

 \begin{proof}[Proof of Theorem \ref{thm-even-assymptotics}] If $|G|=2$, then
$\log_2\log_2 |G|=0$ and $\mathsf s_{\{\pm 1\}}(G)=3=\exp(G)+\log_2|G|$.
 Thus the Theorem holds for any constant $C_1$, and so we may assume $|G|\geq
4$. In this case, $\log_2\log_2|G|\geq 1$, and thus it suffices to
prove the existence of $C_r$ when $\exp(G)\geq n_0$ is sufficiently large, as then
 $\mathsf s_{\{\pm 1\}}(G)\leq \exp(G)+\log_2|G|+C_r\log_2\log_2|G|+C'_r$,
where $C'_r\geq 0$ is the maximum of $\mathsf s_{\{\pm 1\}}(G)$
over all $G$
of rank $\rk(G)=r$ and even exponent
 $\exp(G)<n_0$, and replacing $C_r$ by $C_r+C'_r$ gives the
desired constant that works for all $G$.

The rank $r\geq 1$ will remain fixed throughout the argument.
 Let $G$ be a finite abelian group of rank $r$ and exponent $n=\exp(G)$ even.
Let $m=2^{r+1}$. Note that $m$ depends only on $r$,
 and can thus be treated as a constant with regard to asymptotics. We divide
the proof into four parts.

\medskip

 \textbf{Step 1:} There exists a constant $C>0$, dependent only on $r$, so
that any sequence $S\in \Fc(G)$
 with $|S|\geq Cn^{r/(r+1)}$ contains a $\{\pm 1\}$-weighted zero-sum
subsequence of length $m$.

 First let us see that it suffices to prove that $|S|\geq C'n^{r/(r+1)}$
implies $S$ contains an $\{\pm 1\}$-weighted zero-sum subsequence $T$
of length $|T|\in \{2^1,2^2,\ldots,2^{r+1}\}$.
 Indeed, if we know this to be true, then, for any sequence $S\in \Fc(G)$
with $$|S|\geq (C'+(r+1)2^{r+1})n^{r/(r+1)}\geq C'n^{r/(r+1)}+(r+1)2^{r+1},$$
 we can repeatedly apply this result to $S$ to pull off disjoint weighted
 zero-sum subsequences $T_1,\ldots,T_l$ with $T_1\cdot \ldots \cdot T_l|S$,
\be\label{whooft}|T_1\cdot \ldots T_l|>(r+1)2^{r+1},\ee
 and $|T_i|\in \{2^1,\,2^2,\ldots,2^{r+1}\}$ for all $i$. If $S$ contains no
such subsequence of length $2^{r+1}$, then less than $2^{r+1-j}$ of the $T_i$
 can be of length $2^j$, for $j=1,2,\ldots ,r+1$ (else concatenating a
sufficient number of these $T_i$ would yield a weighted zero-sum of the desired
 length $2^{r+1}$). Consequently, $$|T_1\cdot\ldots\cdot
 T_l|<\Sum{j=1}{r+1}2^{r+1-j}\cdot 2^j=(r+1)2^{r+1},$$ contradicting
 \eqref{whooft}. Thus the step follows with constant $C'+(r+1)2^{r+1}$, and
 we see it suffices to prove  $|S|\geq C'n^{r/(r+1)}$ implies $S$ contains an
$\{\pm 1\}$-weighted zero-sum subsequence $T$
of length $|T|\in \{2^1,2^2,\ldots,2^{r+1}\}$, as claimed.

Let $S=g_1\cdot g_2\cdot \ldots\cdot g_{v}$, where $g_i\in G$. Let
$X$ be the collection of all subsets $I\subseteq [1,v]$ having
cardinality $2^r$. Recall that we associate each $I\subseteq [1,v]$
with the indexed subsequence $S_I:=\prod_{i\in I}g_i$ of $S$. If
$\sigma(S_I)=\sigma(S_J)$ for distinct $I,\,J\in X$, then, by
discarding the commonly indexed terms (as we have done several times
before in Section \ref{sec-small-values}), we obtain a $\{\pm
1\}$-weighted zero-sum subsequence $S_{I\setminus J}\cdot
S_{J\setminus I}$ of $S$ with length $$|I|+|J|-2|I\cap J|=m-2|I\cap
J|.$$ Assuming by contradiction that $S$ contains no $\{\pm
1\}$-weighted zero-sum subsequence $T$ with length $|T|\in
\{2^1,2^2,\ldots,2^{r+1}\}$ and recalling that $m=2^{r+1}$, we
conclude that $$|I\cap J|\in L:=[0,2^r-1]\setminus
\left(\frac{m}{2}-\{2^0,2^1,\ldots,2^{r}\}\right)$$ whenever
$\sigma(S_I)=\sigma(S_J)$ with $I,\,J\in X$ distinct. (Note that
$|L|=2^r-r-1$.) This allows us to give an upper bound on how many
distinct $I\in X$ can have equal corresponding sums. Indeed, Theorem
\ref{Ray-Chaudhuri-Wilson-thm} shows that there can be at most
$\binom{v}{|L|}=\binom{v}{2^r-r-1}$ indexing sets from $X$ having
equal corresponding sums.
 Since $|X|=\binom{v}{2^r}$,
this implies there are at least
$$\binom{v}{2^r}\bigg/\binom{v}{2^r-r-1}\gg v^{r+1}$$ distinct
values attained by the $\sigma(S_I)$ with $I\in X$. Since there are
at most $|G|\leq \exp(G)^r=n^r$ values in total, we conclude that
$$n^r\gg v^{r+1},$$ which implies $v<Cn^{r/(r+1)}$ for some
constant $C>0$ (the above asymptotic notation holds for $v$
sufficiently large with respect to $r$, which is a fixed constant).
Thus, if $|S|\geq Cn^{r/(r+1)}$, then $S$ must contain a weighted
zero-sum subsequence of one of the desired lengths, completing the
step as noted earlier.

\medskip

\textbf{Step 2: } There exists a constant $C'>0$, dependent only on $r$, so that any sequence $S\in \Fc(G)$
with $|S|\geq C'n^{r/(r+1)}$ contains a $\{\pm 1\}$-weighted zero-sum subsequence $T$ of
length $|T|\equiv n\ (\mo\ m)$ and $|T|\leq (r+1)m$.

 The proof is a variation on that of Step 1. Let $S=g_1\cdot g_2\cdot
\ldots\cdot g_{v}$, where $g_i\in G$. Let $\alpha\in [1,m]$
 be the integer such that $n\equiv \alpha\ (\mo\ m)$. Note, since $n$ and $m$ are
both even, that $\alpha$ must be an {\it even} number---this is the only
 place where the hypothesis regarding the parity of $n$ will be used.
 Let $X$ be the collection of all subsets $I\subseteq [1,v]$ having
cardinality $\frac12(\alpha+rm)$, which is an integer as both $\alpha$
 and $m$ are even. Recall that we associate each $I\subseteq [1,v]$ with the
indexed subsequence $S_I:=\prod_{i\in I}g_i$ of $S$.
  If $\sigma(S_I)=\sigma(S_J)$ for distinct $I,\,J\in X$, then, by discarding
the commonly indexed terms, we obtain a $\{\pm 1\}$-weighted zero-sum
  subsequence $S_{I\setminus J}\cdot S_{J\setminus I}$ of $S$ with length
$$|I|+|J|-2|I\cap J|=\alpha+rm-2|I\cap J|.$$
  Assuming by contradiction that $S$ contains no $\{\pm 1\}$-weighted
zero-sum subsequence $T$ with length
   $\{\alpha,\alpha+m,\alpha+2m,\ldots,\alpha+rm\}$ and recalling that
$m=2^{r+1}$, we conclude that
   $$|I\cap J|\in L:=\l[0,\frac{\alpha+rm}2-1\r]\setminus \l\{(r-j)\frac{m}{2}\mid
j=0,1,\ldots,r\r\}$$
 whenever $\sigma(S_I)=\sigma(S_J)$ with $I,\,J\in X$ distinct; note
$|L|=\frac12(\alpha+rm)-r-1$.
 This allows us to give an upper bound on how many distinct $I\in X$ can have
equal corresponding sums. Indeed, Theorem \ref{Ray-Chaudhuri-Wilson-thm}
 shows that there can be at most
 $\binom{v}{|L|}=\binom{v}{(\alpha+rm)/2-r-1}$ indexing sets from $X$
having equal corresponding sums.
Since $|X|=\binom{v}{(\alpha+rm)/2}$, this implies there are at least
  $$\binom{v}{(\alpha+rm)/2}\bigg/\binom{v}{(\alpha+rm)/2-r-1}\gg
v^{r+1}$$ distinct values attained by the $\sigma(S_I)$ with $I\in
X$.
 Since there are at most $|G|\leq \exp(G)^r=n^r$ values in total, we conclude
 that  $$n^r\gg v^{r+1},$$ which implies $v< C'n^{r/(r+1)}$ for some
 constant $C'>0$. Thus, if $|S|\geq C'n^{r/(r+1)}$, then $S$ must contain a
weighted zero-sum subsequence of one of the desired lengths, completing the step.

\medskip

\textbf{Step 3: } There exists a constant $C''>0$, dependent only on $r$, so that
any sequence $S\in \Fc(G)$ with $|S|\geq \log_2|G|+C''\log_2\log_2|G|$
contains a $\{\pm 1\}$-weighted zero-sum subsequence $T$
with length $|T|\equiv 0\ (\mo\ m)$ and $|T|\leq \log_2|G|+C''\log_2\log_2|G|$.


Suppose we can show that, for  any sequence $S\in \Fc(G)$ with \be\label{booboot}|S|
=v\geq \log_2 |G|+C'\log_2\log_2 |G|\ee and \be\label{trikesike}v\equiv 0\ (\mo\ 2^{r+2}=2m),\ee there is
 a weighted zero-sum subsequence of $S$ with length congruent to $0$ modulo $m$. Then, since any sequence $S$
 with $|S|\geq \log_2 |G|+C'\log_2\log_2 |G|+2m$ contains a subsequence $S'$
 with $$\log_2 |G|+C'\log_2\log_2 |G|\leq   |S'|\leq  \log_2 |G|+C'\log_2\log_2 |G|+2m+1$$
 that also satisfies \eqref{trikesike}, and since any weighted zero-sum subsequence of $S'$
 has length trivially bounded from above by $|S'|\leq \log_2 |G|+C'\log_2\log_2 |G|+2m+1$,
 we see that the step holds setting $C''=C'+2m+1$ (in view of $\log_2\log_2|G|\geq 1$).
 We proceed to show this supposition true.

 To that end,
let $S=g_1\cdot g_2\cdot \ldots\cdot g_{v}$, where $g_i\in G$, be a sequence satisfying \eqref{trikesike},
in which case $\lfloor\frac{v+1}{2}\rfloor\equiv 0\ (\mo\ m)$. Let $X$ be the collection of all subsets $I\subseteq [1,v]$
having cardinality $|I|\equiv 0\ (\mo\ m)$.
 In view of \eqref{trikesike} and Lemma \ref{sbinom-lemma}, we find that  \be\label{somebound}|X|
 \geq \frac{2^v}{m}=2^{v-r-1}.\ee
 Recall that we associate each $I\subseteq [1,v]$ with the indexed subsequence $S_I:=\prod_{i\in I}g_i$ of $S$.

If $\sigma(S_I)=\sigma(S_J)$ for distinct $I,\,J\in X$, then, by discarding the commonly indexed terms,
we obtain a $\{\pm 1\}$-weighted
zero-sum subsequence $S_{I\setminus J}\cdot S_{J\setminus I}$ of length $$|I|+|J|-2|I\cap J|.$$
Hence, since $|I|+|J|\equiv 0+0= 0\ (\mo\ m)$,
we see that $S_{I\setminus J}\cdot S_{J\setminus I}$ will be a $\{\pm 1\}$-weighted zero-sum of
length congruent to $0$ modulo $m$
provided $|I\cap J|\equiv 0\ (\mo\ m/2)$. Therefore, assuming to the contrary tat this is not the case, we conclude that
 $|I\cap J|\in L$, where $L=\{1,2,3,\ldots,2^r-1\}+2^r\Z$, whenever $\sigma(S_I)=\sigma(S_J)$ with $I,\,J\in X$ distinct.

This
 allows us to give an upper bound on how many distinct $I\in X$ can have equal corresponding sums.
 Indeed, since all $I\in X$ have $|I|\equiv 0\ (\mo\ m=2^{r+1})$,
  we see that all $I\in X$ have $I\in K+2^r\Z$, where $K=\{0\}$. Moreover, the polynomial $f(x)=\prod_{i=1}^{2^r-1}(x-i)$
  shows that $0$ can be separated from $\{1,2,3,\ldots,2^r-1\}+2^r\Z$ with respect to $p=2$
  using a polynomial of degree $D=2^r-1$.
  Thus, applying Theorem \ref{thm-primepower} with $q=p^k=2^r=m/2$ and using \eqref{somebound},
  we see that there are
  at least $$2^{v-r-1}\bigg/\Sum{i=0}{D}\binom{v}{i}\gg 2^{v-r-1}/v^{2^r-1}$$
distinct values attained by the $\sigma(S_I)$ with $I\in X$. Hence, since there are at most $|G|$ values in total,
 we conclude that \be\label{pieprinkle}|G|\geq C\frac{2^{v-r-1}}{v^{D}}\ee for some $C>0$ when $v\geq v_0$,
 where $v_0>0$ is some constant depending on the fixed constant $r$
(using $D=2^r-1$).

Recall, since $|G|\geq 4$, that \be\label{stuffsee}\log_2 |G|\geq 2\und \log_2 |G|\geq \log_2\log_2|G|\geq 1.\ee
Let $$\gamma=D+\max\{0,\log_2(1/C)\}\in\Z^+.$$ Now $2^{x}$ is larger than $(x+\gamma+r+2)^D$ for sufficiently large $x$.
Thus, let $y\geq v_0$ be an integer such that
\be\label{calcdant}2^{x}>(x+\gamma+r+2)^D\quad\mbox{ for all } x\geq y\ee and consider $C'=y+\gamma+r+1$.

Suppose $|S|=v\geq \log_2 |G|+C'\log_2\log_2 |G|$.
 Then $$v=\log_2|G|+(x+\gamma+r+1)\log_2\log_2 |G|\geq v_0$$ for some real number $x\geq y$,
 and using \eqref{stuffsee}
we derive that  \ber \nn C2^{v-r-1}&=&C|G|2^{(x+\gamma+r+1)\log_2\log_2 |G|-r-1}\geq C|G|2^{(x+\gamma)\log_2\log_2 |G|}
\\ &=&
C|G|(\log_2|G|)^{x+\gamma}
= C|G|(\log_2|G|)^{D}(\log_2|G|)^{x+\max\{0,\log_2(1/C)\}}\nn
\\ &\geq& C|G|(\log_2|G|)^{D}2^{x+\max\{0,\log_2(1/C)\}}
\geq 2^x(\log_2|G|)^{D}|G|\label{calcdant-swo},
\eer and that, again using \eqref{stuffsee},
\ber \nn v^D&=&(\log_2|G|+(x+\gamma+r+1)\log_2\log_2 |G|)^D
 \\&\leq& ((x+\gamma+r+2)\log_2|G|)^D=(x+\gamma+r+2)^D(\log_2|G|)^D.\label{calcdant-swu}
\eer Combining \eqref{calcdant-swo} and \eqref{calcdant-swu} and using \eqref{calcdant} and $x\geq y$, it follows that
 $$C\frac{2^{v-r-1}}{v^D}\geq \frac{2^x(\log_2|G|)^D|G|}{(x+\gamma+r+2)^D(\log_2|G|)^D}
 =\frac{2^x}{(x+\gamma+r+2)^D}|G|>|G|,$$
 contradicting \eqref{pieprinkle}. Thus we see the constant $C'$ for \eqref{booboot} exists,
 completing the step as remarked earlier.

\medskip

\textbf{Step 4: }There exists a constant $C_r>0$, dependent only on $r$, so that, for sufficiently large $n$,
any sequence $S\in \Fc(G)$
 with $|S|\geq n+\log_2|G|+C_r\log_2\log_2|G|$ contains a $\{\pm 1\}$-weighted zero-sum subsequence $T$ with
 of length $|T|=n$.

Note that this step will complete the proof, for as remarked at the beginning of the proof, it suffices to
prove the theorem for sufficiently large $n$.  For this reason, we may also assume $n\geq (r+1)m$.
Let $C$, $C'$ and $C''$ be the respective constants from Steps 1, 2 and 3.
 We will show that $S\in \Fc(G)$ contains a $\{\pm 1\}$-weighted zero-sum subsequence of
 length $n$ provided the length of $S$ satisfies
 the following three bounds:
\ber \label{sojo1}|S|&\geq& C'n^{r/(r+1)},
\\\label{sojo2} |S|&\geq& Cn^{r/(r+1)}+(r+1)m+\log_2|G|+C''\log_2\log_2|G|,\\
\label{sojo3}|S|&\geq& n-m+\log_2|G|+C''\log_2\log_2|G|.\eer
Since, for sufficiently large $n$, the bound given in \eqref{sojo3} is the maximum of the three bounds,
we will subsequently be able to
conclude $|S|\geq n-m+\log_2|G|+C''\log_2\log_2|G|$, for large $n$, implies $S$
contains a weighted zero-sum subsequence of length $n$, completing the proof.
We continue by showing \eqref{sojo1}--\eqref{sojo3} indeed guarantee a length $n$ weighted zero-sum subsequence.

In view of \eqref{sojo1} and Step 2, we see that $S$ contains a weighted zero-sum subsequence $R_0$ with
\be\label{thinklusq}|R_0|\equiv n\,(\mo\ m)\und |R_0|\leq (r+1)m\leq n.\ee In view of \eqref{sojo2},
we see that repeated application of Step 1 to $R_0^{-1}S$ yields series of length $m$ weighted zero-sum subsequences,
 enough so that there exists
a subsequence $R$ of $R_0^{-1}S$ with $$|R|\geq \log_2 |G|+C''\log_2 \log_2|G|\und |R|\equiv 0\ (\mo\ m)$$ such that,
for every $k\in [0,|R|]\cap m\Z$,
there is a $\{\pm 1\}$-weighted zero-sum subsequence $T_k$ of $R$ with length $|T_k|=k$.
 Choose such a subsequence $R$ of $R_0^{-1}S$ with length $|R|$
maximal. Since $|R_0|\equiv n\ (\mo\ m)$ with $|R_0|\leq n$ (in view of \eqref{thinklusq}),
 we see that $n=|R_0|+ym$ for some $y\in \N$.
Thus $|R_0R|\leq n-m$, else the proof is complete.

In view of $|R_0R|\leq n-m$ and \eqref{sojo3}, we see that $$|R_0^{-1}R^{-1}S|\geq \log_2|G|+C''\log_2\log_2|G|.$$
Hence, applying Step 3
to  $R_0^{-1}R^{-1}S$, we find a nontrivial weighted zero-sum subsequence $T$ of $R_0^{-1}R^{-1}S$
with $|T|\equiv 0\ (\mo\ m)$
and $$|T|\leq \log_2|G|+C''\log_2\log_2|G|.$$
We claim that $RT$ contradicts the maximality of $|R|$, which, once shown true,
will provide the concluding contradiction for the proof.

Since $|R|\equiv |T|\equiv 0\ (\mo\ m)$, we have $|RT|\equiv 0\ (\mo\ m)$,
while $$|TR|\geq |R|\geq \log_2 |G|+C''\log_2 \log_2|G|.$$
For every $k\in [0,|R|]\cap m\Z$, the weighted zero-sum subsequence $T_k$ divides $R$, and hence also $RT$,
and is of length $|T_k|=k$.
Since $$|R|\geq \log_2 |G|+C''\log_2 \log_2|G|\geq |T|-m,$$ it follows that every $t\in [|R|+m,|R|+|T|]\cap m\Z$
can be written in the form $t=k+|T|$
with $k\in [|R|-|T|+m,|R|]\cap m\Z\subseteq [0,|R|]\cap m\Z$. Hence the subsequence $T_kT$ of $RT$
is a weighted zero-sum subsequence of length
$t\in [|R|+m,|R|+|T|]\cap m\Z$, which shows that the subsequence  $RT$ of $R_0^{-1}S$ indeed contradicts the
maximality of $|R|$, completing the proof. \qedsymbol
\end{proof}

\medskip

\noindent{\bf Acknowledgment}. The authors would like to thank Dr.
Hao Pan at Nanjing University for his helpful comments on the proof
of Theorem 1.1.

\end{document}